\documentclass[11pt]{article}
\newcommand{\qed}{{\hfill\rule{4pt}{7pt}}}
\usepackage{amsmath, epsfig, cite}
\usepackage{amssymb}
\usepackage{amsfonts}
\usepackage{latexsym}
\newtheorem{thm}{Theorem}

\newtheorem{lem}[thm]{Lemma}

\def\Z{\mathbb{Z}}
\def\pf{\noindent {\it Proof.} }

\begin{document}

\begin{center}
{\Large\bf The Number of Convex Polyominoes and the Generating
Function of Jacobi Polynomials}
\end{center}

\vskip 2mm
\centerline{Victor J. W. Guo$^1$ and Jiang Zeng$^{1,2}$}

\begin{center} $^1$ Center for Combinatorics, LPMC\\
Nankai University, Tianjin 300071, People's Republic of China\\
{\tt jwguo@eyou.com}
\end{center}

\begin{center}
$^2$ Institut Girard Desargues,
Universit\'e Claude Bernard (Lyon I)\\
F-69622 Villeurbanne Cedex, France \\
{\tt zeng@desargues.univ-lyon1.fr}
\end{center}

\vskip 0.7cm \noindent{\bf Abstract.} {\small Lin and Chang gave a
generating function of convex polyominoes with an $m+1$ by $n+1$
minimal bounding rectangle.  Gessel  showed that their result
implies that the number of such polyominoes is
$$
\frac{m+n+mn}{m+n}{2m+2n\choose 2m}-\frac{2mn}{m+n}{m+n\choose m}^2.
$$
 We show that this result can be derived from some binomial
 coefficients identities related to the generating function of Jacobi polynomials.
}

\begin{quote}
Some (binomial coefficients) identities arise from
alternative solutions of combinatorial problems and incidentally
give added significance to doing problems the ``hard" way.
\flushright{--- {\sc J. Riordan}}
\end{quote}



\noindent{\bf Abbreviated title}: The Number of Convex Polyominoes

\section{Introduction}
A {\em polyomino}  is a connected union of squares in the plane
whose vertices are lattice points. A polyomino is called {\em
convex} if its intersection with any horizontal or vertical line
is either empty or a line segment. Any convex polyomino has a
minimal bounding rectangle whose perimeter is the same as that of
the polyomino. Delest and Viennot~\cite{DV} found a generating
function for counting convex polyominoes by perimeter and showed
that the number of convex polyominoes with perimeter $2n+8$, for
$n\geq 0$, is
\begin{equation}\label{eq:DV}
(2n+11)4^n-4(2n+1){2n\choose n}.
\end{equation}
Later, Lin and Chang~\cite{LC} gave a generating function for the
number of convex polyominoes with an $(m+1)\times (n+1)$ minimal
bounding rectangle, and Gessel~\cite{Gessel} showed that their
result implies that the number of such polyominoes is
\begin{equation}\label{eq:Ge}
\frac{m+n+mn}{m+n}{2m+2n\choose 2m}-\frac{2mn}{m+n}{m+n\choose m}^2,
\end{equation}
which is easily seen to give a refinement of Delest and Viennot's
formula.

Since Gessel~\cite{Gessel} (see also Bousquet-M\'elou~\cite{BM})
derived  \eqref{eq:Ge} from the  generating function of Lin and
Chang~\cite{LC} (see also  Bousquet-M\'elou and
Guttman~\cite{BMG}), it would be interesting to find an
independent proof of \eqref{eq:Ge}. For Delest and Viennot's
formula~\eqref{eq:DV} such a proof  was already given by
Kim~\cite{Kim}. The aim of this paper is to provide such a proof
for \eqref{eq:Ge} by generalizing Kim's elementary approach. It
turns out that the resulting binomial coefficients identities are
related to the generating function of Jacobi polynomials.

In the next section, we translate the enumeration of convex
polyominoes with fixed minimal bounding rectangle as that of two
pairs of non intersecting lattice paths, which results to evaluate
a quadruple sum of binomial coefficients. In Section 3, we
establish some binomial coefficients identities which lead to the
evaluation of the desired sums.

\section{Non intersecting lattice paths and determinant formula}
\begin{figure}[h]
\setlength{\unitlength}{0.5cm}
\begin{picture}(12.5,10)(-8,0)
\put(-1,8.5){$n+1$}\put(13,0.5){$m+1$}
 \qbezier[40](0,.5)(6,.5)(12.5,.5)
\qbezier[40](0,8)(6,8)(12.5,8) \qbezier[40](0,0.5)(0,4)(0,8)
\qbezier[40](12.5,0.5)(12.5,4)(12.5,8)
 \put(0,3.5){\circle*{.3}}
\put(0,3.5){\line(0,1){2.5}}\put(0,6){\line(1,0){3}}
\put(4,6){$L_4$} \put(3,6){\line(0,1){1}}\put(3,7){\line(1,0){1}}
\put(4,7){\line(0,1){1}}\put(4,8){\line(1,0){3}}
\put(7,8){\circle*{.3}}
\put(0,3.5){\line(1,0){2.5}}\put(2.5,3.5){\line(0,-1){2}} \put(3,
2 ){$L_1$} \put(7.5,0.5){\circle*{.3}}
\put(2.5,1.5){\line(1,0){4}}\put(6.5,1.5){\line(0,-1){1}}
\put(6.5,0.5){\line(1,0){4}}
\put(10.5,0.5){\line(0,1){2}} \put(10.5,2.5){\line(1,0){2}}
\put(11,3){$L_3$} \put(12.5,2.5){\line(0,1){2}}
\put(12.5,4.5){\circle*{.3}}
\put(7,8){\line(0,-1){1.5}}\put(7,6.5){\line(1,0){3}}
\put(9,5){$L_2$}
\put(10,6.5){\line(0,-1){2}}\put(10,4.5){\line(1,0){2.5}}
\put(12.5,4.5){\circle*{.3}}
\put(6.5,-.5){$a_1$}\put(-1,3.5){$b_1$}
\put(7,8.5){$m+1-a_2$}\put(13,4.5){$n+1-b_2$}
\put(6.5,0.5){\circle*{.3}}\put(0,4.5){\circle*{.3}}
\put(6,8){\circle*{.3}}\put(12.5,3.5){\circle*{.3}}
\end{picture}
\vskip 2mm
\caption{A convex polyomino with an $m+1$ by $n+1$ minimal
bounding rectangle.\label{fig:convex}}
\end{figure}
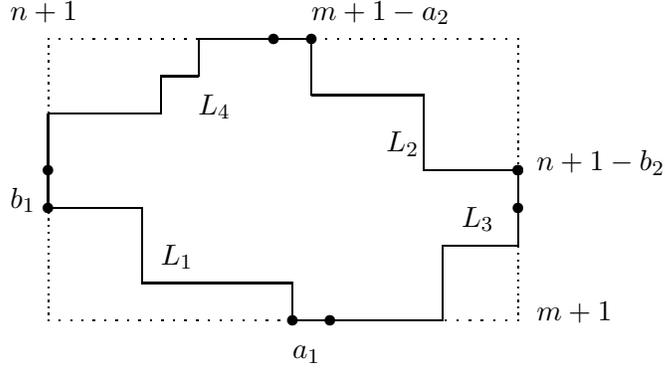

A {\em lattice path} is a sequence of points $(s_0,s_1,\ldots,
s_n)$ in the plan $\Z^2$ such that either $s_i-s_{i-1}=(1,0),\,
(0,1)$ for all $i=1,\ldots, n$ or $s_i-s_{i-1}=(1,0),\, (0,-1)$
for all $i=1,\ldots, n$. Let ${\cal P}_{m,n}$ be the set of convex
polyominoes with an $m+1$ by $n+1$ minimal bounding rectangle. As
illustrated in Figure~\ref{fig:convex}, any polyomino in ${\cal
P}_{m,n}$ can be characterized by 4 lattice paths $L_1$, $L_2$,
$L_3$ and $L_4$ which are given by
\begin{align*}
L_1&\colon\ (0,b_1)\longrightarrow (a_1,0),\\
L_2&\colon\ (m+1-a_2,n+1)\longrightarrow (m+1, n+1-b_2),\\
L_3&\colon\ (a_1+1,0)\longrightarrow (m+1,n-b_2),\\
L_4&\colon\ (0,b_1+1)\longrightarrow (m-a_2,n+1).
\end{align*}

Note that  a polyomino in ${\cal P}_{m,n}$ is convex if and only
if the two lattice paths $L_1,L_2$ (resp. $L_3,L_4$) don't
intersect. The following lemma can be readily proved by switching
the {\em tails} of two lattice paths, which is also a special case
of a more general result~\cite{GV}.
\begin{lem}
Let $a,b,c$ and $d$ be non negative integers such that $a'>a$,
$b>b'$ , $c>a$, $d>b$, $c'>a'$ and $d'>b'$. Then the number of
pairs of non intersecting lattice paths $({\cal P}_1, {\cal P}_2)$
such that  ${\cal P}_1\colon (a,b)\longrightarrow (c,d)$ and
${\cal P}_2\colon (a',b')\longrightarrow (c',d')$ is given by
\[
{c-a+d-b\choose c-a}{c'-a'+d'-b'\choose c'-a'}
-{c-a'+d-b'\choose c-a'}{c'-a+d'-b\choose c'-a}.
\]
\end{lem}

It follows that the cardinality of ${\cal P}_{m,n}$ is given by
\begin{align}
&\sum_{a_1,a_2=0}^m\sum_{b_1,b_2=0}^n
 \left[{a_1+b_1-2\choose a_1-1}{a_2+b_2-2\choose a_2-1}
-{a_1+a_2+n-m-2\choose n-1}{b_1+b_2+m-n-2\choose m-1} \right] \nonumber\\[5pt]
&\cdot
 \left[{m+n-a_2-b_1\choose m-a_2}{m+n-a_1-b_2\choose m-a_1}
 -{m+n-a_1-a_2\choose n+1}{m+n-b_1-b_2\choose  m+1}\right],\label{eq:initial}
\end{align}
Note that in \eqref{eq:initial} we have adopted the
convention that ${-2\choose -1}=1$, which corresponds to
$a_1=b_1=0$ or $a_2=b_2=0$. In this case the path $L_1$ or $L_2$
is a point.

We next split the sum in \eqref{eq:initial} into three terms: the
$a_1=a_2=b_1=b_2=0$ term,
\begin{equation}\label{eq:S0}
S_0 = {m+n\choose m}^2-{m+n\choose m-1}{m+n\choose n-1},
\end{equation}
the $a_1=b_1=0$ or $a_2=b_2=0$  terms,
\begin{align*}
S_1&= 2\sum_{a=1}^m\sum_{b=1}^n {a+b-2\choose a-1}
 \left[{m+n-a\choose m-a}{m+n-b\choose m}-
{m+n-a\choose n+1}{m+n-b\choose m+1}\right],
\end{align*}
and the sum in \eqref{eq:initial} for $a_i$ and $b_i\geq 1$. This
third term can be split into four more terms obtained from the
product:
$(a_{11}a_{22}-a_{12}a_{21})(b_{11}b_{22}-b_{12}b_{21})=S_2-S_3-S_4+S_5$.
The last sum $S_5$ is 0 because the sum of the numerator
parameters of binomials coefficients  are less than that of the
denominator parameters.

 We now proceed to evaluate or simplify $S_1$, $S_2$, $S_3$ and $S_4$
 using  the Chu-Vandermonde formula:
$$ {}_2F_1\left(\begin{array}{cc}
  -n, & a \\
 & \hspace{-15pt} c
\end{array}; 1\right):=\sum_{k\geqslant0}\frac{(-n)_k(a)_k}{(c)_k
k!}=\frac{(c-a)_n}{(c)_n},
$$
where $(x)_n=x(x+1)\ldots (x+n-1)$ for $n\geq 1$ and $(x)_0=1$.

\begin{itemize}
\item Applying  the Chu-Vandermonde formula to the $b$-sums for
$S_1$ yields
\[
S_1=2\sum_{a=1}^m \left[{m+n-a\choose n}{m+n+a-1\choose n-1}-
{m+n-a\choose n+1}{m+n+a-1\choose n-2}\right].
\]
As
\begin{align*}
&{m+n-a\choose n}{m+n+a-1\choose n-1}
-{m+n-a\choose n+1}{m+n+a-1\choose n-2} \\[5pt]
&={m+n-a+1\choose n+1}{m+n+a-1\choose n-1}
 -{m+n-a\choose n+1}{m+n+a\choose n-1},
\end{align*}
by telescoping it follows that
\begin{equation}
S_1=2{m+n\choose n+1}{m+n\choose n-1}.
\end{equation}

\item Consider now the second sum $S_2$:
$$
S_2=\sum_{a_1,a_2=1}^m\sum_{b_1,b_2=1}^n{a_1+b_1-2\choose a_1-1}
{a_2+b_2-2\choose a_2-1}{m+n-a_2-b_1\choose m-a_2}
{m+n-a_1-b_2\choose m-a_1}.
$$
By the Chu-Vandermonde formula we have
\begin{align*}
\sum_{b_1=1}^n{a_1+b_1-2\choose a_1-1}{m+n-a_2-b_1\choose m-a_2}
&={m+n+a_1-a_2-1\choose n-1}, \\[5pt]
\sum_{b_2=1}^n{a_2+b_2-2\choose a_2-1}{m+n-a_1-b_2\choose m-a_1}
&={m+n-a_1+a_2-1\choose n-1}.
\end{align*}
Hence
$$
S_2=\sum_{a_1,a_2=1}^m{m+n+a_1-a_2-1\choose n-1}{m+n-a_1+a_2-1\choose n-1}.
$$
Setting $a=a_1-a_2$ we can rewrite the above sum as
\begin{align*}
S_2&= \sum_{a=1-m}^{m-1}\#\{(a_1,a_2)\in [1,m]^2 \mid a_1-a_2=a\}
{m+n+a-1\choose n-1}{m+n-a-1\choose n-1}\\[5pt]
&= \sum_{a=-m}^m(m-|a|){m+n+a-1\choose n-1}{m+n-a-1\choose n-1}\\[5pt]
&= m\sum_{a=-m}^m{m+n+a-1\choose n-1}{m+n-a-1\choose n-1}\\[5pt]
&\quad\ -2\sum_{a=1}^m a{m+n+a-1\choose n-1}{m+n-a-1\choose n-1}.
\end{align*}

By the Chu-Vandermonde formula  we have
\begin{align*}
m\sum_{a=-m}^{m}{m+n+a-1\choose n-1}{m+n-a-1\choose n-1}
=m{2m+2n-1\choose 2n-1}.
\end{align*}
Since
\begin{align*}
&2a{m+n+a-1\choose n-1}{m+n-a-1\choose n-1} \\[5pt]
&=n{m+n+a-1\choose n}{m+n-a\choose n}
-n{m+n+a\choose n}{m+n-a-1\choose n},
\end{align*}
telescoping yields
\[
\sum_{a=1}^{m}2a{m+n+a-1\choose n-1}{m+n-a-1\choose n-1}
=n{m+n\choose n}{m+n-1\choose n}.
\]
Hence
\begin{equation}\label{eq:S2}
S_2=\frac{mn}{m+n}{2m+2n\choose 2m}-\frac{mn}{m+n}{m+n\choose
m}^2.
\end{equation}

\item Look at the term $S_3$:
\begin{align*}
S_3&=\sum_{a_1,a_2=1}^m\sum_{b_1,b_2=1}^n{a_1+a_2+n-m-2\choose n-1}
{b_1+b_2+n-m-2\choose m-1}\\[5pt]
&\hspace{3cm}\cdot {m+n-a_2-b_1\choose m-a_2} {m+n-a_1-b_2\choose m-a_1}.
\end{align*}
Summing the $a_2$-sum and $b_2$-sum by the Chu-Vandermonde formula
yields
\begin{align*}
\sum_{a_2=1}^{m}{a_1+a_2+n-m-2\choose n-1}{m+n-a_2-b_1\choose m-a_2}
&={2n+a_1-b_1-1\choose a_1-1},\\[5pt]
\sum_{b_2=1}^{n}{b_1+b_2+n-m-2\choose m-1}{m+n-a_1-b_2\choose m-a_1}
&={2m-a_1+b_1-1\choose b_1-1}.
\end{align*}
Hence, replacing $a_1$ and $b_1$ by $a$ and $b$ respectively we
get
\begin{align*}
S_3&=\sum_{a=1}^m\sum_{b=1}^n{2m-a+b-1\choose b-1}{2n+a-b-1\choose a-1}\\[5pt]
&=\sum_{a=1}^m\sum_{b=1}^n{m+n+a-b-1\choose m+a-1}{m+n-a+b-1\choose n+b-1},
\end{align*}
by the substitutions $a\leftarrow m-a+1$ and $b\leftarrow n-b+1$.
\item Finally we have
$$
S_4=\sum_{a_1,a_2=1}^m\sum_{b_1,b_2=1}^n{a_1+b_1-2\choose a_1-1}
{a_2+b_2-2\choose a_2-1}{m+n-a_1-a_2\choose n+1}
{m+n-b_1-b_2\choose m+1}.
$$
Summing the $a_1$-sum and $b_2$-sum by the Chu-Vandermonde formula
yields
\begin{align*}
\sum_{a_1=1}^{m}{a_1+b_1-2\choose a_1-1}{m+n-a_1-a_2\choose n+1}
&={m+n-a_2+b_1-1\choose n+b_1+1},\\[5pt]
\sum_{b_2=1}^{n}{a_2+b_2-2\choose a_2-1}{m+n-b_1-b_2\choose m+1}
&={m+n+a_2-b_1-1\choose m+a_2+1}.
\end{align*}
Substituting $a_2$ and $b_1$ by $a$ and $b$ we obtain
$$
S_4=\sum_{a=1}^{m-2}\sum_{b=1}^{n-2}{m+n+a-b-1\choose m+a+1}
{m+n-a+b-1\choose n+b+1},
$$
for the summand is zero if $a=m-1, m$ or $b=n-1,n$.
\end{itemize}
We shall evaluate $S_3$ and $S_4$ in the next section.
\section{Jacobi polynomials and evaluation of $S_3$ and $S_4$}
Set
\[
\Delta:=\sqrt{1-2x-2y-2xy+x^2+y^2}.
\]
The following identity is equivalent to the generating function of
Jacobi polynomials:
\begin{align}
\sum_{m,n=0}^{\infty} {m+n+\alpha\choose m}{m+n+\beta\choose n}x^m y^n
&=\frac{2^{\alpha+\beta}}{\Delta(1-x+y+\Delta)^\alpha(1+x-y+\Delta)^\beta}.
\label{eq:pqdelta}
\end{align}
The reader is referred to \cite[p.~298]{AAR} and~\cite[p.~271]{Ra} for
two classical analytical proofs and to \cite{FL} for  a
combinatorial proof.

Applying the operator $x\frac{\partial}{\partial x}+y\frac{\partial}{\partial y}+2$
to the $\alpha=\beta=1$ case of \eqref{eq:pqdelta} yields:
\begin{equation}\label{eq:delta}
\sum_{m,n\geq 1}\frac{m+n}{2}{m+n-1\choose m}{m+n-1\choose
n}x^my^n=\frac{xy}{\Delta^3}.
\end{equation}

\begin{thm}\label{thm:mnabsqare}
We have
\begin{equation}\label{eq:S3}
S_3=\frac{mn}{2(m+n)}{m+n\choose m}^2.
\end{equation}
\end{thm}
\pf  Consider the generating function of $S_4$:
\begin{align*}
F(x,y)&:=\sum_{m,n=0}^{\infty}\sum_{a=1}^{m}\sum_{b=1}^{n}
{m+n-a+b-1\choose m-a}{m+n+a-b-1\choose n-b}x^my^n \\[5pt]
&=\sum_{a=1}^{\infty}\sum_{b=1}^{\infty}x^ay^b\sum_{m=a}^{\infty}\sum_{n=b}^{\infty}
{m+n-a+b-1\choose m-a}{m+n+a-b-1\choose n-b}x^{m-a} y^{n-b} \nonumber\\[5pt]
&=\sum_{a,b=1}^{\infty}x^ay^b\sum_{m,n=0}^{\infty}
{m+n+2b-1\choose m}{m+n+2a-1\choose n}x^{m} y^{n}.
\end{align*}
Applying \eqref{eq:pqdelta} to the inner double sum yields
\[
F(x,y)=\sum_{a,b=1}^{\infty}x^ay^b
\frac{2^{2a+2b-2}}{\Delta(1-x+y+\Delta)^{2b-1}(1+x-y+\Delta)^{2a-1}}=\frac{xy}{\Delta^3}.
\]
The theorem follows then from \eqref{eq:delta}. \qed

\begin{thm}There holds
\begin{equation}\label{eq:S4}
S_4={m+n\choose m}^2+{m+n\choose m-1}{m+n\choose
n-1}+\frac{mn}{2(m+n)}{m+n\choose m}^2-{2m+2n\choose 2n}.
\end{equation}
\end{thm}
\pf  Consider the generating function of $S_4$:
\begin{align*}
G(x,y)&:=\sum_{m,n=0}^{\infty}\sum_{a=1}^{m-2}\sum_{b=1}^{n-2}
{m+n-a+b-1\choose m-a-2}{m+n+a-b-1\choose n-b-2}x^my^n \\[5pt]
&=\sum_{a=1}^{\infty}\sum_{b=1}^{\infty}
\sum_{m=a+2}^{\infty}\sum_{n=b+2}^{\infty}
{m+n-a+b-1\choose m-a-2}{m+n+a-b-1\choose n-b-2}x^{m} y^{n} \nonumber\\[5pt]
&=\sum_{a,b=1}^{\infty}x^{a+2}y^{b+2} \sum_{m,n=0}^{\infty}
{m+n+2b+3\choose m}{m+n+2a+3\choose n}x^{m} y^{n}.
\end{align*}
Applying \eqref{eq:pqdelta} to the inner double sum yields
\begin{align*}
G(x,y)&=\sum_{a,b=1}^{\infty}x^{a+2}y^{b+2}
\frac{2^{2a+2b+6}}{\Delta(1-x+y+\Delta)^{2b+3}(1+x-y+\Delta)^{2a+3}}\\[5pt]
&=\frac{16x^3y^3}{\Delta^3(1-x-y+\Delta)^4}.
\end{align*}
Set
\[
f(x,y):=\sum_{m,n=0}^{\infty}{m+n\choose m}x^my^n=\frac{1}{1-x-y}.
\]
By bisecting twice, we get the terms of even powers of $x$ and $y$
in $f(x,y)$:
\begin{align*}
\sum_{m,n=0}^{\infty}{2m+2n\choose 2m}x^{2m}y^{2n}
=\frac{1}{4}(f(x,y)+f(-x,y)+f(x,-y)+f(-x,-y))
\end{align*}
i.e.,
\[
\sum_{m,n=0}^{\infty}{2m+2n\choose
2m}x^{m}y^{n}=\frac{1-x-y}{\Delta^2}.
\]
Now, the $\alpha=\beta=0$ and $\alpha=\beta=2$ cases of \eqref{eq:pqdelta} read:
\begin{align*}
\sum_{m,n=0}^\infty{m+n\choose m}^2x^my^n&=\frac{1}{\Delta},\\
\sum_{m,n=1}^\infty{m+n\choose m-1}{m+n\choose n-1}
x^my^n&=\frac{4xy}{\Delta(1-x-y+\Delta)^2}.
\end{align*}
As
\[
\frac{16x^3y^3}{\Delta^3(1-x-y+\Delta)^4}
=\frac{1}{\Delta}+\frac{4xy}{\Delta(1-x-y+\Delta)^2}+\frac{xy}{\Delta^3}
-\frac{1-x-y}{\Delta^2},
\]
extracting the coefficients of $x^my^n$ in the above equation
completes the proof. \qed

Summarizing, formula \eqref{eq:Ge} follows then from
\eqref{eq:S0}--\eqref{eq:S2}, \eqref{eq:S3} and \eqref{eq:S4}.

\end{document}